\documentclass[12pt]{article}
\usepackage{latexsym, amssymb}
\usepackage{mathbbol}

\DeclareMathAlphabet\gothic{U}{euf}{m}{n}

\makeatletter
\def\eqnarray{\stepcounter{equation}\let\@currentlabel=\theequation
\global\@eqnswtrue
\tabskip\@centering\let\\=\@eqncr
$$\halign to \displaywidth\bgroup\hfil\global\@eqcnt\z@
  $\displaystyle\tabskip\z@{##}$&\global\@eqcnt\@ne
  \hfil$\displaystyle{{}##{}}$\hfil
  &\global\@eqcnt\tw@ $\displaystyle{##}$\hfil
  \tabskip\@centering&\llap{##}\tabskip\z@\cr}

\def\endeqnarray{\@@eqncr\egroup
      \global\advance\c@equation\m@ne$$\global\@ignoretrue}

\def\@yeqncr{\@ifnextchar [{\@xeqncr}{\@xeqncr[5pt]}}
\makeatother

\begin{document}
\bibliographystyle{tom}

\newtheorem{lemma}{Lemma}[section]
\newtheorem{thm}[lemma]{Theorem}
\newtheorem{cor}[lemma]{Corollary}
\newtheorem{voorb}[lemma]{Example}
\newtheorem{rem}[lemma]{Remark}
\newtheorem{prop}[lemma]{Proposition}
\newtheorem{stat}[lemma]{{\hspace{-5pt}}}
\newtheorem{obs}[lemma]{Observation}
\newtheorem{defin}[lemma]{Definition}

\newenvironment{remarkn}{\begin{rem} \rm}{\end{rem}}
\newenvironment{exam}{\begin{voorb} \rm}{\end{voorb}}
\newenvironment{defn}{\begin{defin} \rm}{\end{defin}}
\newenvironment{obsn}{\begin{obs} \rm}{\end{obs}}

\newcommand{\gota}{\gothic{a}}
\newcommand{\gotb}{\gothic{b}}
\newcommand{\gotc}{\gothic{c}}
\newcommand{\gote}{\gothic{e}}
\newcommand{\gotf}{\gothic{f}}
\newcommand{\gotg}{\gothic{g}}
\newcommand{\gothh}{\gothic{h}}
\newcommand{\gotk}{\gothic{k}}
\newcommand{\gotm}{\gothic{m}}
\newcommand{\gotn}{\gothic{n}}
\newcommand{\gotp}{\gothic{p}}
\newcommand{\gotq}{\gothic{q}}
\newcommand{\gotr}{\gothic{r}}
\newcommand{\gots}{\gothic{s}}
\newcommand{\gotu}{\gothic{u}}
\newcommand{\gotv}{\gothic{v}}
\newcommand{\gotw}{\gothic{w}}
\newcommand{\gotz}{\gothic{z}}
\newcommand{\gotA}{\gothic{A}}
\newcommand{\gotB}{\gothic{B}}
\newcommand{\gotG}{\gothic{G}}
\newcommand{\gotL}{\gothic{L}}
\newcommand{\gotS}{\gothic{S}}
\newcommand{\gotT}{\gothic{T}}

\newcounter{teller}
\renewcommand{\theteller}{\Roman{teller}}
\newenvironment{tabel}{\begin{list}%
{\rm \bf \Roman{teller}.\hfill}{\usecounter{teller} \leftmargin=1.1cm
\labelwidth=1.1cm \labelsep=0cm \parsep=0cm}
                      }{\end{list}}

\newcounter{tellerr}
\renewcommand{\thetellerr}{(\roman{tellerr})}
\newenvironment{subtabel}{\begin{list}%
{\rm  (\roman{tellerr})\hfill}{\usecounter{tellerr} \leftmargin=1.1cm
\labelwidth=1.1cm \labelsep=0cm \parsep=0cm}
                         }{\end{list}}
\newenvironment{ssubtabel}{\begin{list}%
{\rm  (\roman{tellerr})\hfill}{\usecounter{tellerr} \leftmargin=1.1cm
\labelwidth=1.1cm \labelsep=0cm \parsep=0cm \topsep=1.5mm}
                         }{\end{list}}

\newcommand{\Ni}{{\bf N}}
\newcommand{\Ri}{{\bf R}}
\newcommand{\Ci}{{\bf C}}
\newcommand{\Ti}{{\bf T}}
\newcommand{\Zi}{{\bf Z}}
\newcommand{\Fi}{{\bf F}}

\newcommand{\proof}{\mbox{\bf Proof} \hspace{5pt}} 
\newcommand{\remark}{\mbox{\bf Remark} \hspace{5pt}}
\newcommand{\ruimte}{\vskip10.0pt plus 4.0pt minus 6.0pt}

\newcommand{\simh}{{\stackrel{{\rm cap}}{\sim}}}
\newcommand{\ad}{{\mathop{\rm ad}}}
\newcommand{\Ad}{{\mathop{\rm Ad}}}
\newcommand{\Aut}{\mathop{\rm Aut}}
\newcommand{\arccot}{\mathop{\rm arccot}}
\newcommand{\capp}{{\mathop{\rm cap}}}
\newcommand{\rcapp}{{\mathop{\rm rcap}}}
\newcommand{\diam}{\mathop{\rm diam}}
\newcommand{\divv}{\mathop{\rm div}}
\newcommand{\codim}{\mathop{\rm codim}}
\newcommand{\RRe}{\mathop{\rm Re}}
\newcommand{\IIm}{\mathop{\rm Im}}
\newcommand{\Tr}{{\mathop{\rm Tr}}}
\newcommand{\Vol}{{\mathop{\rm Vol}}}
\newcommand{\card}{{\mathop{\rm card}}}
\newcommand{\supp}{\mathop{\rm supp}}
\newcommand{\sgn}{\mathop{\rm sgn}}
\newcommand{\essinf}{\mathop{\rm ess\,inf}}
\newcommand{\esssup}{\mathop{\rm ess\,sup}}
\newcommand{\Int}{\mathop{\rm Int}}
\newcommand{\Leibniz}{\mathop{\rm Leibniz}}
\newcommand{\lcm}{\mathop{\rm lcm}}
\newcommand{\loc}{{\rm loc}}

\newcommand{\mod}{\mathop{\rm mod}}
\newcommand{\spann}{\mathop{\rm span}}
\newcommand{\one}{\mathbb{1}}

\hyphenation{groups}
\hyphenation{unitary}

\newcommand{\tfrac}[2]{{\textstyle \frac{#1}{#2}}}

\newcommand{\cb}{{\cal B}}
\newcommand{\cc}{{\cal C}}
\newcommand{\cd}{{\cal D}}
\newcommand{\ce}{{\cal E}}
\newcommand{\cf}{{\cal F}}
\newcommand{\ch}{{\cal H}}
\newcommand{\ci}{{\cal I}}
\newcommand{\ck}{{\cal K}}
\newcommand{\cl}{{\cal L}}
\newcommand{\cm}{{\cal M}}
\newcommand{\co}{{\cal O}}
\newcommand{\cs}{{\cal S}}
\newcommand{\ct}{{\cal T}}
\newcommand{\cx}{{\cal X}}
\newcommand{\cy}{{\cal Y}}
\newcommand{\cz}{{\cal Z}}

\newcommand{\wtozp}{W^{1,2}\raisebox{10pt}[0pt][0pt]{\makebox[0pt]{\hspace{-34pt}$\scriptstyle\circ$}}}
\newlength{\hightcharacter}
\newlength{\widthcharacter}
\newcommand{\covsup}[1]{\settowidth{\widthcharacter}{$#1$}\addtolength{\widthcharacter}{-0.15em}\settoheight{\hightcharacter}{$#1$}\addtolength{\hightcharacter}{0.1ex}#1\raisebox{\hightcharacter}[0pt][0pt]{\makebox[0pt]{\hspace{-\widthcharacter}$\scriptstyle\circ$}}}
\newcommand{\cov}[1]{\settowidth{\widthcharacter}{$#1$}\addtolength{\widthcharacter}{-0.15em}\settoheight{\hightcharacter}{$#1$}\addtolength{\hightcharacter}{0.1ex}#1\raisebox{\hightcharacter}{\makebox[0pt]{\hspace{-\widthcharacter}$\scriptstyle\circ$}}}
\newcommand{\scov}[1]{\settowidth{\widthcharacter}{$#1$}\addtolength{\widthcharacter}{-0.15em}\settoheight{\hightcharacter}{$#1$}\addtolength{\hightcharacter}{0.1ex}#1\raisebox{0.7\hightcharacter}{\makebox[0pt]{\hspace{-\widthcharacter}$\scriptstyle\circ$}}}

\thispagestyle{empty}

\begin{center}
{\Large\bf Flows and invariance for elliptic operators} \\[5mm]
\large A.F.M. ter Elst$^1$, Derek W. Robinson$^2$ and Adam Sikora$^3$

\end{center}

\vspace{5mm}

\begin{center}
{\bf Abstract}
\end{center}

\begin{list}{}{\leftmargin=1.8cm \rightmargin=1.8cm \listparindent=10mm 
   \parsep=0pt}
\item
Let $S$ be the submarkovian semigroup on $L_2(\Ri^d)$ generated by a
self-adjoint, second-order, divergence-form, elliptic operator $H$ 
with $W^{1,\infty}$ coefficients $c_{kl}$. 
Further let $\Omega$ be an open subset of $\Ri^d$. 
Under mild conditions we prove that $S$ leaves $L_2(\Omega)$
invariant if, and only if, it is invariant under the flows 
generated by the vector fields $\sum_{l=1}^d c_{kl} \, \partial_l$
for all $k$.
\end{list}

\vspace{7cm}
\noindent
March 2009

\vspace{5mm}

\vspace{5mm}

\noindent
{\bf Home institutions:}    \\[3mm]
\begin{tabular}{@{}cl@{\hspace{10mm}}cl}
1. & Department of Mathematics  & 
  2. & Centre for Mathematics   \\
& University of Auckland   & 
  & \hspace{15mm} and its Applications  \\
& Private bag 92019 & 
  & Mathematical Sciences Institute  \\
& Auckland & 
  & Australian National University  \\
& New Zealand  & 
  & Canberra, ACT 0200  \\
& & 
  & Australia \\[8mm]
  3. & Department of Mathematics& 
 {}\\
& Macquarie University & 
  {}\\
&Sydney, NSW 2109 & 
 {}\\
& Australia & 
 {} \\
& {} & 
 {} \\
& {} & 
  & {} 

\end{tabular}

\newpage
\setcounter{page}{1}

\section{Introduction}\label{Sflow1}

Let $S$ be a submarkovian semigroup on $L_2(\Ri^d)$ generated by a 
self-adjoint second-order elliptic operator $H$ in divergence form. 
If the operator is strongly elliptic  then $S$ acts ergodically, 
i.e.\ there are no non-trivial $S$-invariant subspaces of $L_2(\Ri^d)$.
Nevertheless there are many examples of degenerate elliptic operators 
for which there are 
subspaces $L_2(\Omega)$  invariant under the action of $S$ (see, for 
example, \cite{ERSZ1} \cite{ERSZ2}
\cite{RSi} \cite{ER29}).
Our aim is to examine operators with coefficients  which are  
Lipschitz continuous and  characterize the  $S$-invariance of 
$L_2(\Omega)$ by the invariance under a family of associated  flows. 
In order to formulate our main result we need  
some further notation.

First  define the positive symmetric operator $H_0$ with
domain $D(H_0)=C_c^\infty(\Ri^d)$ and action
\[
H_0\varphi= -\sum^d_{k,l=1} \partial_k \, c_{kl} \, \partial_l\varphi
\]
where the coefficients $c_{kl}=c_{lk}\in W^{1,\infty}(\Ri^d)$ are 
real and  $C=(c_{kl})$
is a  positive-definite matrix over $\Ri^d$.
Then the corresponding quadratic form $h_0$ given by 
\[
h_0(\varphi)
= \sum^d_{k,l=1}(\partial_k\varphi, c_{kl} \, \partial_l\varphi)
\]
with domain $ D(h_0)=C_c^\infty(\Ri^d)$ is closable.
The closure $h=\overline{h_0}$ determines  in a canonical manner a 
positive self-adjoint extension $H$ of $H_0$, the Friedrichs' extension \cite{Friedr2}
(see, for example, \cite{RN}, \S 124, or \cite{Kat1}, Chapter~VI).
The closed form   $h$ is  a Dirichlet form and the self-adjoint 
semigroup $S$ generated by $H$ is automatically submarkovian
(for details on Dirichlet forms and submarkovian semigroups see 
\cite{FOT} or \cite{BH}).
We call $H$ the degenerate elliptic operator with coefficients $(c_{kl})$.

Secondly, if $b_1,\ldots,b_d \in W^{1,\infty}(\Ri^d)$ then 
the first-order partial differential operator 
\[
\varphi \mapsto \sum^d_{k=1} b_k \, \partial_k \varphi
   - \frac{1}{2} \sum^d_{k=1}(\partial_k b_k) \, I
\]
with domain $C_c^\infty(\Ri^d)$ is 
essentially skew-adjoint (see, for example, \cite{Rob7}, Theorem~3.1).
Therefore the principal part is closable and 
generates a  positive, continuous,  one-parameter group 
on $L_2(\Ri^d)$.
We refer to such a group as flows.
Specifically we are interested in the flows associated with the coefficients $(c_{kl})$ of $H$.
For all $k \in \{ 1,\dots, d \} $ let $Y_k$ 
denote the $L_2$-closures of the first-order
partial differential operator
\[
\varphi \mapsto \sum^d_{l=1} c_{kl} \, \partial_l \varphi
\]
with domain $C_c^\infty(\Ri^d)$.
Then  denote by $T^{(k)}$ the flows generated by the $Y_k$.
The operators $Y_k$  were used by  Ole\u{\i}nik and Radkevi\v{c} \cite{OR} to analyze   
hypoellipticity and subellipticity properties of degenerate elliptic operators 
$H$ with $C^\infty$-coefficients $c_{kl}$ (see \cite{JSC} for a review of these and related results).
We, however, use the flows to characterize the invariant subspaces of the semigroup generated by $H$.

\begin{thm}\label{tsep1} 
Let $\Omega$ be a measurable subset of $\Ri^d$.
Consider the following conditions.
\begin{tabel}
\item\label{tsep1-1}
$S_t L_2(\Omega) \subseteq L_2(\Omega)$ for all $t>0$.
\item\label{tsep1-2}
$T^{(k)}_t L_2(\Omega) = L_2(\Omega)$ for all $k \in \{ 1,\ldots,d \} $ and $t\in\Ri$.
\end{tabel}
Then {\rm\ref{tsep1-1}}$\Rightarrow${\rm\ref{tsep1-2}}.
Moreover,  
if 
$C_c^\infty(\Ri^d)$ is a core for $H$, or, if
 $\Omega$ is open  and the boundary $\partial \Omega$ of $\Omega$
is  $($locally$)$ Lipschitz then {\rm\ref{tsep1-1}}$\Leftrightarrow${\rm\ref{tsep1-2}}.
\end{thm}

Recall that the   open set $\Omega$ is defined to have a 
(locally) Lipschitz boundary if for 
every $y\in \partial \Omega$ there exist an isometry $\Psi \colon \Ri^d\to \Ri^d$,
a real function $\tau \in W^{1,\infty}( \Ri^{d-1} )$
and an $r > 0$ such that 
\begin{equation}
 \Omega\cap B_y(r) = \{ \Psi(x_1,x') :  (x_1,x') \in \Ri \times \Ri^{d-1}, \; 
                               \tau(x')  < x_1 \} \cap B_y(r)
\label{eSflow1;1}
\end{equation}
where $B_y(r)=\{x\in\Ri^d: \|x-y\|<r\}$.
Thus in a neighbourhood of $y$ the boundary $\partial\Omega$ 
of $\Omega$ is the graph of a Lipschitz function $\tau$, up to an 
isometry $\Psi$.

There are two variations of the theorem which will be established in the course of its proof.

First, for all $\psi \in C_c^\infty(\Ri^d)$ define $Y_\psi$ 
as the $L_2$-closure of the first-order partial differential operator
\[
\varphi \mapsto \sum^d_{k,l=1} (\partial_k \psi) \, c_{kl} \, \partial_l \varphi
\]
with domain $C_c^\infty(\Ri^d)$ and let $T^\psi$ be the associated flow.
Then invariance of $L_2(\Omega)$ under the $T^{(k)}$ is equivalent to invariance under
the family of flows $T^\psi$.
More precisely one has the following.

\begin{prop}\label{psep1} 
Let $\Omega$ be a measurable subset of $\Ri^d$.
The following conditions are equivalent:
\begin{tabel}
\item\label{psep1-4}
$T^\psi_t L_2(\Omega) = L_2(\Omega)$ for all $\psi \in C_c^\infty(\Ri^d)$ and $t\in\Ri$.
\item\label{psep1-2}
$T^{(k)}_t L_2(\Omega) = L_2(\Omega)$ for all $k \in \{ 1,\ldots,d \} $ and $t\in\Ri$.
\end{tabel}
\end{prop}
This will be established in Section~\ref{Sflow2}.

Secondly, the condition that $C_c^\infty(\Ri^d)$ is a core for $H$ does not follow in general from 
the assumption that the coefficients are in  $ W^{1,\infty}(\Ri^d)$.
The one-dimensional example considered in \cite{ERSZ1}, Section~5 
gives a counterexample.
Specifically, let $\delta \in [1/2,\infty \rangle$ and
$H=-d\,c\,d$ with $c(x)= |x|^{2\delta}(1+x^2)^{-\delta}$.
Then $c\in W^{1,\infty}(\Ri)$ but $C_c^\infty(\Ri)$ is a
core of $H$ if and only if  $\delta\geq 3/4$
by the arguments in \cite{CMP}, Proposition~3.5.
(See also \cite{RSi3}.)
In particular it is not a core if $\delta\in[1/2,3/4\rangle$. 
Nevertheless it follows that $C_c^\infty(\Ri^d)$ is a core for $H$ if  
$c_{kl}\in W^{2,\infty}(\Ri^d)$ (see,  \cite{Rob7} Section~6, or \cite{ER27}
Proposition~2.3).
Moreover, the core condition can be derived  from weaker smoothness 
assumptions on the  $c_{kl}$ (see Section~\ref{Sflow4}).

\section{Flows} \label{Sflow2} 

In this section we derive some properties of the flows defined in 
Section~\ref{Sflow1}  and prove Proposition~\ref{psep1}.
Although we deal primarily with the flows on $L_2(\Ri^d)$ we will need,
in Section~\ref{S3}
some properties of their extensions to $L_\infty(\Ri^d)$.
Therefore we  begin by summarizing some general features of the flows.

Let $b_1,\ldots,b_d \in W^{1,\infty}(\Ri^d)$ and define 
$Y$  as  the $L_2$-closure of the first-order differential operator
$\varphi \mapsto \sum_{k=1}^d b_k \, \partial_k\varphi$ and domain 
$W^{1,2}(\Ri^d)$.
Further let $T$ denote the flow generated by $Y$.
Then for all $p\in[1,\infty]$ the group $T$ leaves the subspace
$L_2(\Ri^d)\cap L_p(\Ri^d)$ of $L_2(\Ri^d)$ invariant and 
$T$ extends from $L_2(\Ri^d)\cap L_p(\Ri^d)$  to a flow 
$T^{[p]}$ on $L_p(\Ri^d)$ such that $T^{[p]}$ is strongly continuous if 
$p\in[1,\infty\rangle$ and $T^{[\infty]}$ is weakly$^*$ continuous.
The groups act in a consistent and compatible manner on the $L_p$-spaces.
Moreover, $T^{[\infty]}$ is a group of automorphisms of  $L_\infty(\Ri^d)$, i.e.\
$T^{[\infty]}_t(\psi\, \varphi)= (T^{[\infty]}_t \psi) \, (T^{[\infty]}_t \varphi)$
for all $\psi,\varphi\in L_\infty(\Ri^d)$ and $t \in \Ri$.
Then since the $L_\infty$-functions are multipliers on the $L_p$-spaces
one deduces that 
\begin{equation}
T^{[p]}_t(\tau \, \varphi)
= (T^{[\infty]}_t \tau) \, (T^{[p]}_t\varphi)
\label{eSflow2;5}
\end{equation}
for all $\tau \in L_\infty(\Ri^d)$, $\varphi\in L_p(\Ri^d)$, $p\in[1,\infty]$ 
and $t \in \Ri$.
If $Y_{[p]}$ is the generator of $T^{[p]}$ then $W^{1,p}(\Ri^d) \subset D(Y_{[p]})$
and $Y_{[p]} \varphi = \sum_{k=1}^d b_k \, \partial_k\varphi$  
for all $\varphi \in W^{1,p}(\Ri^d)$.

These properties  depend critically on the fact that  $Y$  is a first-order partial 
differential operator with coefficients $b_k\in W^{1,\infty}(\Ri^d)$.
They can be verified either by general arguments of functional analysis 
(see, for example,  \cite{Robm}, Theorem~V.4.1) or by methods of ordinary differential
equations.
The crucial observation in the latter context is that if 
$\varphi\in C_c^\infty(\Ri^d)$ then 
$(T_t\varphi)(x)=\varphi(\omega_t(x))$ where 
$t \mapsto \omega_t(x)$ is the unique solution of  the differential equation
$(d/dt)\omega_t(x)=b(\omega_t(x))$, with initial value $\omega_0(x)=x$ (see, for example, \cite{Hil}, Chapters~2 and 3).

Our first result is an approximation result which will be needed on $L_2(\Ri^d)$ but whose proof extends to the $L_p$-spaces.

\begin{prop} \label{pord901}
Let $p \in [1,\infty]$.
Let $Y_{[p]}$ denote the generator of the flow $T^{[p]}$ on $L_p(\Ri^d)$.
Further let $\tau \in C_c^\infty(\Ri^d)$ with $\int \tau = 1$ and 
for all $n \in \Ni$ define $\tau_n \in C_c^\infty(\Ri^d)$ by 
$\tau_n(x) = n^d \, \tau(n \, x)$.

Then $\lim_{n \to \infty} Y_{[p]}(\tau_n * \varphi) = Y_{[p]} \varphi$
in $L_p(\Ri^d)$ for all $\varphi \in D(Y_{[p]})$ if $p < \infty$.
If $p = \infty$ then 
$\lim_{n \to \infty} Y_{[\infty]}(\tau_n * \varphi) = Y_{[\infty]} \varphi$ weakly$^*$
in $L_\infty(\Ri^d)$ for all $\varphi \in D(Y_{[\infty]})$
\end{prop}
\proof\
First, for all $n \in \Ni$ define the bounded operator
$B_n \colon L_p \to L_p$ by
\[
B_n \varphi
= \sum_{k=1}^d \tau_n * ( (\partial_k \, b_k) \varphi)
   + \sum_{k=1}^d \int dy \, (\partial_k \tau_n)(y) \, 
        \Big( (I - L_y) b_k \Big) \, (L_y \varphi) 
\;\;\; ,  \]
where $L$ denotes  the left regular representation of $\Ri^d$, i.e.\
$(L_y \psi)(x) = \psi(x-y)$.
Secondly, if $\varphi \in C_c^\infty$ and $n \in \Ni$ then
\begin{eqnarray*}
Y_{[p]}(\tau_n * \varphi)
& = & \sum_{k=1}^d b_k \int dy \, \tau_n(y) \, L_y \partial_k \, \varphi  \\
& = & \sum_{k=1}^d \int dy \, \tau_n(y) \, (b_k - L_y b_k) \, L_y \partial_k \, \varphi
  + \sum_{k=1}^d b_k \int dy \, \tau_n(y) \, L_y (b_k \, \partial_k \, \varphi)  
\;\;\; .  
\end{eqnarray*}
The second term equals $\tau_n * Y_{[p]} \varphi$.
For the first term use
$L_y \partial_k \varphi = - \frac{\partial}{\partial y_k} \, L_y \varphi$.
Therefore  integration by parts gives
\[
Y_{[p]}(\tau_n * \varphi) - \tau_n * Y_{[p]} \varphi
= \sum_{k=1}^d \int dy \, \frac{\partial}{\partial y_k}
    \Big( \tau_n(y) \, (b_k - L_y b_k) \Big) (L_y \varphi)  
= B_n \varphi
\;\;\; .  \]
Since $B_n$ is bounded one deduces by density that 
\begin{equation}
Y_{[p]}(\tau_n * \varphi) - \tau_n * Y_{[p]} \varphi = B_n \varphi
\label{epord901;1}
\end{equation}
for all $n \in \Ni$ and $\varphi \in D(Y_{[p]})$.

Thirdly, it follows from the definition of $B_n$ that 
\begin{eqnarray*}
\|B_n \varphi\|_p
& \leq & \sum_{k=1}^d \Big( \|(\partial_k \, b_k) \, \varphi\|_p
   + \int dy \, |(\partial_k \tau_n)(y)| \, 
        \|\Big( (I - L_y) b_k \Big) \, (L_y \varphi)\|_p \Big)  \\
& \leq & \sum_{k=1}^d \|b_k\|_{W^{1,\infty}} \, \|\varphi\|_p
   + \sum_{k=1}^d \int dy \, |(\partial_k \tau_n)(y)| \, 
        \|(I - L_y) b_k\|_\infty \, \|\varphi\|_p 
\end{eqnarray*}
for all $n \in \Ni$ and $\varphi \in L_p$.
But $\|(I - L_y) b_k\|_\infty \leq |y| \, \|b_k\|_{W^{1,\infty}}$
and $\int dy \, |(\partial_k \tau_n)(y)| \, |y| = \int dy \, |(\partial_k \tau)(y)| \, |y|$.
Therefore $\|B_n \varphi\|_p \leq M \, \|\varphi\|_p$ uniformly 
for all $n \in \Ni$ and $\varphi \in L_p$, where 
$M = \sum_{k=1}^d (1 + \int dy \, |(\partial_k \tau)(y)| \, |y|) \, \|b_k\|_{W^{1,\infty}}$.
The conclusion holds for all $p\in[1,\infty]$.
So $B_1, B_2,\ldots$ are equicontinuous.

Next assume $p < \infty$.
If $\varphi \in W^{1,p}$ then 
$\lim_{n \to \infty} \tau_n * \varphi = \varphi$
in $W^{1,p}$.
Consequently, $\lim_{n \to \infty} Y_{[p]}\,(\tau_n * \varphi) = Y_{[p]}\, \varphi$ strongly in $L_p$.
Moreover, $\lim_{n \to \infty} \tau_n * (Y_{[p]} \,\varphi) = Y_{[p]}\, \varphi$ strongly in $L_p$.
Therefore $\lim_{n \to \infty} B_n \varphi = 0$ in $L_p$ for all 
$\varphi \in W^{1,p}$ by (\ref{epord901;1}).
Since $W^{1,p}$ is strongly dense in $L_p$ and $B_1,B_2,\ldots$ are equicontinuous
it follows that $\lim_{n \to \infty} B_n \varphi = 0$ in $L_p$ for all 
$\varphi \in L_p$.
Finally, let $\varphi \in D(Y_{[p]})$.
Then one establishes from (\ref{epord901;1}) that 
$\lim_{n \to \infty} Y_{[p]}\,(\tau_n * \varphi) 
= \lim_{n \to \infty} (\tau_n * Y_{[p]}\, \varphi + B_n \varphi)
= Y_{[p]} \,\varphi$ in $L_p$.

The argument for $p=\infty$ is very similar.
If $\varphi\in W^{1,\infty}$ then  $\lim \tau_n * \varphi = \varphi$
and $\lim \partial_k \tau_n * \varphi = \partial_k\varphi$ weakly$^*$.
Therefore 
$\lim Y_{[\infty]}\,(\tau_n * \varphi) = Y_{[\infty]}\, \varphi$ weak$^*$ on $L_\infty$.
Then since $W^{1,\infty}$ is weakly$^*$  dense in $L_\infty$ and 
$B_1,B_2,\ldots$ are equicontinuous the desired conclusion  follows 
as before.\hfill$\Box$

\ruimte

Now we return to consideration of the vector fields $Y_1,\ldots,Y_d$
defined in Section~\ref{Sflow1} acting on $L_2(\Ri^d)$.

\begin{cor} \label{cflow310}
Let  $\tau$ and $\tau_n$ be as in Proposition~{\rm \ref{pord901}}.
Then for all $\varphi \in \bigcap_{k=1}^d D(Y_k)$ one has
$\lim_{n \to \infty} Y_k(\tau_n * \varphi) = Y_k \varphi$
for all $k \in \{ 1,\ldots,d \} $.
\end{cor}

Note that convolution with $\tau_n$ maps $L_2(\Ri^d)$ into $W^{\infty,2}(\Ri^d)$ 
so the corollary
establishes that $W^{\infty,2}(\Ri^d)$ is  a simultaneous core for  
the $Y_1,\ldots,Y_d$.

Now we turn to the proof of Proposition~\ref{psep1}.
Note  that 
if $T$ is a flow  with generator $Y$ then $T$-invariance of $L_2(\Omega)$  is equivalent to the 
the commutation of $Y$ and the operator of  multiplication with $\one_\Omega$,
i.e.\  if  $\varphi \in D(Y)$ then $\one_\Omega \varphi\in  D(Y)$ and 
$Y(\one_\Omega \, \varphi) = \one_\Omega \, Y \varphi$.

\ruimte

\noindent
{\bf Proof of Proposition~\ref{psep1}\hspace{5pt} }
``\ref{psep1-4}$\Rightarrow$\ref{psep1-2}''.
Let $k \in \{ 1,\ldots,d \} $ and $U \subset \Ri^d$ a bounded open subset.
There exist $\chi,\psi \in C_c^\infty(\Ri^d)$ such that 
$\chi|_U = \one$ and $\psi(x) = x_k$ for all $x \in \supp \chi$.
Then $Y_k(\chi \varphi) = Y_\psi(\chi \varphi)$ 
for all $\varphi \in C_c^\infty(\Ri^d)$.
Since $\varphi \mapsto \chi \varphi$ is continuous on $D(Y_k)$
and on $D(Y_\psi)$, with the graph norm, it follows from Proposition~\ref{pord901}
that $\chi \varphi \in D(Y_k)$ for all $\varphi \in D(Y_\psi)$.
In particular, if $\varphi \in C_c^\infty(\Ri^d)$ with $\supp \varphi \subset U$
then $\one_\Omega \varphi \in D(Y_\psi)$ and therefore
$\one_\Omega \varphi = \chi \one_\Omega \varphi \in D(Y_k)$.
Moreover, 
$Y_k(\one_\Omega \varphi) 
= Y_\psi(\chi \one_\Omega \varphi)
= \one_\Omega Y_\psi(\chi \varphi)
= \one_\Omega Y_k \varphi$.
Then it follows by continuity that $\one_\Omega \varphi\in  D(Y_k)$ and 
$Y_k(\one_\Omega \, \varphi) = \one_\Omega \, Y_k \varphi$ for all $\varphi \in D(Y_k)$.
Therefore Condition~\ref{psep1-2} is valid.

``\ref{psep1-2}$\Rightarrow$\ref{psep1-4}''.
It follows from Condition~\ref{psep1-2} that 
$\one_\Omega \varphi \in D(Y_k)$ and 
$Y_k (\one_\Omega \varphi) = \one_\Omega Y_k \varphi$
for all $\varphi \in D(Y_k)$.
Let $\psi \in C_c^\infty(\Ri^d)$.
Then $Y_\psi \varphi = \sum_{k=1}^d (\partial_k \psi) \, Y_k \varphi$
for all $\varphi \in C_c^\infty(\Ri^d)$. 
Since the coefficients $c_{kl}$ are in $W^{1,\infty}(\Ri^d)$ it 
follows from Corollary~\ref{cflow310} that 
$\varphi \in D(Y_\psi)$ and 
$Y_\psi \varphi = \sum_{k=1}^d (\partial_k \psi) \, Y_k \varphi$
for all $\varphi \in \bigcap_{k=1}^d D(Y_k)$.
Hence if $\varphi \in C_c^\infty(\Ri^d)$ then 
$\one_\Omega \varphi \in D(Y_\psi)$ and 
$Y_\psi (\one_\Omega \varphi) = \one_\Omega Y_\psi \varphi$.
By density the latter extends to all $\varphi \in D(Y_\psi)$ and 
therefore Condition~\ref{psep1-4} is valid.
\hfill$\Box$

\ruimte

Finally we note that the flows $T^\psi$ can be defined for all $\psi \in W^{2,\infty}(\Ri^d)$
and the conditions of  Proposition~\ref{psep1} are equivalent to invariance of $L_2(\Omega)$
for all  $T^\psi_t$ with $\psi \in W^{2,\infty}(\Ri^d)$ and $t>0$.
This follows from the arguments of the  foregoing proof.

\section{Semigroup invariance}\label{S3}

In this section we prove Theorem~\ref{tsep1}.
First, however, we  observe that Condition~\ref{tsep1-2} of the theorem, the invariance of $L_2(\Omega)$
under the flows $T^{(k)}$  is equivalent to  $T^\psi$-invariance of $L_2(\Omega)$ for all $\psi\in C_c^\infty(\Ri^d)$.
This is a direct consequence of  Proposition~\ref{psep1}  which was established in the previous section.
Therefore in the subsequent discussion we will consider the $T^\psi$-invariance condition.

\ruimte

\noindent
{\bf Proof of Theorem~\ref{tsep1}\hspace{5pt} }
``\ref{tsep1-1}$\Rightarrow$\ref{tsep1-2}''.
It suffices, by the foregoing observation,  to prove the $T^\psi$-invariance of $L_2(\Omega)$ for all $\psi\in C_c^\infty(\Ri^d)$.

First, it follows from the density of $C_c^\infty(\Ri^d)$ in $D(h)$ 
that there exists a unique bilinear map
$\Gamma \colon D(h) \times D(h) \to L_1$, the {\it carr\'e du champ}, such that 
\[
\Gamma(\psi,\varphi)
= \sum_{k,l=1}^d c_{kl} \, (\partial_k \psi) \, (\partial_l \varphi)
\]
for all $\psi,\varphi \in W^{1,2}(\Ri^d)$.
Then $\|\Gamma(\psi,\varphi)\|_1 \leq h(\psi)^{1/2} \, h(\varphi)^{1/2}$
for all $\psi,\varphi \in D(h)$ by the Cauchy--Schwarz inequality.
Moreover,
\begin{equation}
\int  \tau \, \Gamma(\psi,\varphi)
= \frac{1}{2} \Big( h(\tau\psi,\varphi) + h(\psi,\tau\varphi) - h(\tau,\psi\varphi)\Big)
\label{eSord9;1}
\end{equation}
for all $\tau,\psi,\varphi \in C_c^\infty(\Ri^d)$. 
But 
(\ref{eSord9;1}) then extends to all $\tau,\psi,\varphi \in D(h) \cap L_\infty$ by density.

Secondly, the form $h$ is local in the sense that 
$h(\psi,\varphi) = 0$ for all $\psi,\varphi \in D(h)$
with $\psi \, \varphi = 0$
(see \cite{Schm}).
Therefore it follows from (\ref{eSord9;1})  that $\Gamma$ is local in the same sense.

Thirdly, since $L_2(\Omega)$ is  $S$-invariant the operation of 
multiplication by $\one_\Omega$ maps 
$D(h)$ into itself.
Therefore if $\psi,\varphi,\tau \in D(h) \cap L_\infty$
then $\one_\Omega \varphi, \one_\Omega \tau \in D(h) \cap L_\infty$.
By locality of $h$
one deduces from (\ref{eSord9;1}) that 
\begin{eqnarray*}
\int  \tau \, \Gamma(\psi,\one_\Omega \varphi)
& = & \frac{1}{2} \Big( h(\tau\psi,\one_\Omega \varphi) 
   + h(\psi,\tau\one_\Omega \varphi) - h(\tau,\psi \one_\Omega \varphi) \Big)  \\
& = & \frac{1}{2} \Big( h(\one_\Omega \tau\psi,\varphi) 
   + h(\psi,\one_\Omega \tau\varphi) - h(\one_\Omega \tau,\psi\varphi) \Big)  
= \int \one_\Omega  \tau \, \Gamma(\psi,\varphi)
\;\;\; .  
\end{eqnarray*}
Hence $\Gamma(\psi,\one_\Omega \varphi) = \one_\Omega \Gamma(\psi,\varphi)$.
But $D(h) \cap L_\infty$ is dense in $D(h)$.
Therefore $\Gamma(\psi,\one_\Omega \varphi) = \one_\Omega \Gamma(\psi,\varphi)$
for all $\psi,\varphi \in D(h)$.

Now fix $\psi \in C_c^\infty(\Ri^d)$.
Let $\tau \in C_c^\infty(\Ri^d)$.
Then 
\[
((Y_\psi)^* \tau, \eta)
= (\tau, Y_\psi \eta)
= (\tau, \Gamma(\psi, \eta))
\]
for all $\eta \in C_c^\infty(\Ri^d)$.
Since $C_c^\infty(\Ri^d)$ is dense in $D(h)$ one deduces that 
$((Y_\psi)^* \tau, \eta) = (\tau, \Gamma(\psi, \eta))$
for all $\eta \in D(h)$.
Choosing $\eta = \one_\Omega \varphi$ it follows that 
\[
((Y_\psi)^* \tau, \one_\Omega \varphi)
= (\tau, \Gamma(\psi,\one_\Omega \varphi))
= (\one_\Omega \tau, \Gamma(\psi,\varphi))
= (\one_\Omega \tau, Y_\psi \varphi)
= (\tau, \one_\Omega Y_\psi \varphi)
\;\;\; .  \]
Since $C_c^\infty(\Ri^d)$ is a core for $(Y_\psi)^*$ one deduces that 
$\one_\Omega \varphi \in D(Y_\psi)$ and 
$Y_\psi (\one_\Omega \varphi) = \one_\Omega Y_\psi \varphi$.
This conclusion then extends to all $\varphi \in D(Y_\psi)$ by density.
Therefore $L_2(\Omega)$ is invariant under $T^\psi$.

The converse implication \ref{tsep1-2}$\Rightarrow$\ref{tsep1-1}
consists of  two special cases.

{\em Case 1}.  $C_c^\infty(\Ri^d)$ is a core for $H$.

Condition~\ref{tsep1-2} is equivalent to 
$T^\psi$ invariance of $L_2(\Omega)$ for all $\psi\in C_c^\infty(\Ri^d)$ by Proposition~\ref{psep1}.
 Therefore  we assume the latter condition.

Let $\psi,\tau \in C_c^\infty(\Ri^d)$.
Then 
\[
(H \psi, \tau \, \varphi)
= h(\psi, \tau \, \varphi)
= \int \Gamma(\psi, \tau \, \varphi)  
= \int \tau \, \Gamma(\psi,\varphi) + \varphi \, \Gamma(\psi,\tau)
= (\tau,Y_\psi \varphi) + (\varphi, Y_\psi \tau)
\]
for all $\varphi \in C_c^\infty(\Ri^d)$.
Since $C_c^\infty(\Ri^d)$ is dense in $D(Y_\psi)$ one deduces that 
\begin{equation}
(H \psi, \tau \, \varphi)
= (\tau, Y_\psi \varphi) + (\varphi, Y_\psi \tau)
\label{etftos302;1}
\end{equation}
for all $\varphi \in D(Y_\psi)$.

Now let $\psi,\tau,\varphi \in C_c^\infty(\Ri^d)$.
Then by $T^\psi$-invariance of $L_2(\Omega)$ and (\ref{etftos302;1}) one deduces that 
$\one_\Omega \, \varphi \in D(Y_\psi)$ and 
\begin{eqnarray*}
(H \psi, \tau \, \one_\Omega \, \varphi)
& = & (\tau, Y_\psi (\one_\Omega \, \varphi)) 
         + (\one_\Omega \, \varphi, Y_\psi \tau)  \\
& = & (\one_\Omega \, \tau, Y_\psi \varphi) 
      + (\one_\Omega \, \varphi, Y_\psi \tau)  
 =  (\one_\Omega \, \tau, \Gamma(\psi, \varphi)) 
      + (\one_\Omega \, \varphi, \Gamma(\psi, \tau)) 
\;\;\; .
\end{eqnarray*}
Therefore 
\[
|(H \psi, \tau \, \one_\Omega \, \varphi)|
\leq \|\one_\Omega \, \tau\|_\infty \, \|\Gamma(\psi, \varphi)\|_1 
    + \|\one_\Omega \, \varphi\|_\infty \, \|\Gamma(\psi, \tau)\|_1
\leq c \, h(\psi)^{1/2} 
\leq c \, \|(I + H)^{1/2} \psi\|_2
\]
where $c = \|\tau\|_\infty \, h(\varphi)^{1/2} + \|\varphi\|_\infty \, h(\tau)^{1/2}$.
This estimate is uniform for all $\psi \in C_c^\infty(\Ri^d)$.
Since by assumption the space $C_c^\infty(\Ri^d)$ is a core for $D(H)$ it follows
that $\one_\Omega \, \tau \, \varphi \in D(H^{1/2}) = D(h)$ for all 
$\tau,\varphi \in C_c^\infty(\Ri^d)$.
But $\spann (C_c^\infty(\Ri^d) \cdot C_c^\infty(\Ri^d))$ is dense in $D(h)$.
Therefore it follows from \cite{ER29}, Proposition~2.1 III$\Rightarrow$I, that 
$S$ leaves $L_2(\Omega)$ invariant.
This completes the proof of the first case  in the proof of  \ref{tsep1-2}$\Rightarrow$\ref{tsep1-1}.

\medskip

{\em Case 2}.  $\partial \Omega$ is (locally) Lipschitz.

Let $P_\Omega$ be the orthogonal projection of $L_2(\Ri^d)$ onto
$L_2(\Omega)$.
By assumption $T^\psi$ leaves $L_2(\Omega)$ invariant for all 
$\psi\in C_c^\infty(\Ri^d)$.
Hence 
\begin{equation}
T^\psi_t \, P_\Omega = P_\Omega \, T^\psi_t \, P_\Omega
\label{eSflow2;1}
\end{equation}
for all $t \in \Ri$.
Let $B$ denote multiplication by the bounded function
$\sum_{k,l=1}^d (\partial_k \psi) (\partial_l c_{kl}) $ and set 
$M_t=e^{-tB}$ for $\in\Ri$.
Clearly each $M_t$ leaves $L_2(\Omega)$ invariant.
Therefore $(T^\psi_{-t/n} \, M_{-t/n})^n$ leaves $L_2(\Omega)$ 
invariant for all $t \in \Ri$ and $n \in \Ni$.
But $(Y_\psi)^* = - Y_\psi - B$.
Then  the Trotter product formula establishes that 
$(T^\psi_t)^*$  is the strong limit of $(T^\psi_{-t/n} \, M_{-t/n})^n$   
as $n \to \infty$.
So $(T^\psi_t)^*$ leaves $L_2(\Omega)$ invariant.
Hence  $(T^\psi_t)^* \, P_\Omega = P_\Omega \, (T^\psi_t)^* \, P_\Omega$
for all $t \in \Ri$.
Therefore $P_\Omega \, T^\psi_t = P_\Omega \, T^\psi_t \, P_\Omega$
and by (\ref{eSflow2;1}) it follows that 
$T^\psi_t \, P_\Omega = P_\Omega \, T^\psi_t$ for all $t \in \Ri$.
Then 
\[
\one_\Omega \, T^\psi_t \varphi
= P_\Omega \, T^\psi_t \varphi
= T^\psi_t \, P_\Omega \varphi
= T^\psi_t(\one_\Omega \, \varphi)
= (T^{\psi,\infty}_t \one_\Omega) \, (T^\psi_t \varphi)
\]
for all $\varphi \in C_c^\infty(\Ri^d)$ and $t \in \Ri$
where $T^{\psi,\infty}$ denotes the extension of the flow $T^\psi$ to $L_\infty(\Ri^d)$ 
(see Section~\ref{Sflow2}) and we have used (\ref{eSflow2;5}).
Since $T^\psi_t(C_c^\infty(\Ri^d))$ is dense in $L_2(\Ri^d)$
one deduces that $T^{\psi,\infty}_t \one_\Omega = \one_\Omega$ for all 
$t \in \Ri$.

Next let $\varphi \in C_c^\infty(\Ri^d)$.
Then $(Y_\psi)^* \varphi \in L_1(\Ri^d) \cap L_2(\Ri^d)$,
so $(Y^{(\infty)}_\psi)^* \varphi = (Y_\psi)^* \varphi$, 
where $Y^{(\infty)}_\psi$ is the generator of $T^{\psi,\infty}$.
Since 
$((T^{\psi,\infty}_t)^* \varphi, \one_\Omega)
= (\varphi, T^{\psi,\infty}_t \one_\Omega)
= (\varphi, \one_\Omega)$
for all $t \in \Ri$ it follows by differentiation that 
$((Y_\psi)^* \varphi, \one_\Omega) = 0$.
Therefore setting $\Phi_k = \sum_{l=1}^d c_{kl} \, \partial_l \psi$
for $k \in \{ 1,\ldots,d \} $ one has 
\begin{equation}
 \int_\Omega \divv(\varphi \, \Phi)=((Y_\psi)^* \varphi, \one_\Omega) = 0
\;\; \;.
\label{eflux340}
\end{equation}
At this point we use the (local) Lipschitz continuity of $\partial\Omega$.

The Gauss--Green theorem is valid for open sets $\Omega$  with a 
(locally) Lipschitz boundary 
(see, for example, \cite{EvG} page 209).
It states that 
\[
\int_\Omega \divv \Psi
=\int_{\partial \Omega}dS\,  \langle n, \Psi\rangle
\]
for all $\Psi\in W^{1,\infty}(\Ri^d)$ with compact support where 
$\langle\,\cdot\,,\,\cdot\,\rangle$ denotes the inner product on $\Ri^d$,
$dS$ is the  Euclidean measure on $\partial\Omega$ 
and  $n$ is  the   unit outward normal to $\partial\Omega$.
The normal is defined $dS$-almost everywhere.
Thus if one sets $\Psi=\varphi\,\Phi$ with $\varphi\in C_c^\infty(\Ri^d)$
one has
 \[
\int_\Omega \divv(\varphi \, \Phi)
=\int_{\partial \Omega}dS\, \varphi \, \langle n, \Phi\rangle=0
\]
where the last equality uses (\ref{eflux340}).
Since this is valid for all  $\varphi\in C_c^\infty(\Ri^d)$ it follows that 
$\langle n, \Phi \rangle = 0$ almost everywhere on $\partial\Omega$.
Therefore  $\langle (\nabla \psi)(x), C(x) \, n_x \rangle = 0$ 
for  almost every $x \in \partial \Omega$.
But  this is also valid for all $\psi \in C_c^\infty(\Ri^d)$.
Hence one must have 
$C(x) \, n_x = 0$ for almost every $x \in \partial \Omega$.
This corresponds to the condition of zero flux across the boundary as defined 
in \cite{RSi} and then the 
$S$-invariance of $L_2(\Omega)$ follows from Theorem~1.2 of this 
reference.\hfill$\Box$

\ruimte

The  argument  in \cite{RSi} that zero flux implies invariance is somewhat indirect as it first proves that the capacity of $\partial\Omega$
with respect to $h$ is zero and then uses this to deduce the $S$-invariance of $L_2(\Omega)$.
Nevertheless, the same reasoning can be adapted to give a direct proof of the invariance
since the proof can be reduced to a local estimate
as in \cite{RSi}. 
(The latter proof and this proof are an adaption of the argument used to prove Proposition~6.5
in \cite{ERSZ1}.)

First, it suffices to prove that if  $\varphi\in C_c^\infty(\Ri^d)$ then $\one_\Omega\varphi\in D(h)$.
This is a consequence of \cite{ER29} Proposition~2.1 and locality of~$h$.
But this is obvious if the support of $\varphi$ and the boundary are disjoint.
Therefore it suffices to consider $\varphi$ with support close to the boundary $\partial\Omega$.
Then, however, one can use a decomposition of the identity  to reduce to 
the case  $\supp\varphi \subset B_y(r)$ with $y \in \partial\Omega$ and $r > 0$ small.

Secondly, let $\tau$, $\Psi$ be as in (\ref{eSflow1;1}).
Without loss of generality we may assume that $\Psi(x) = x$ for all $x \in \Ri^d$.
For all $n \in \Ni$ define $\psi_n \colon \Ri^d \to \Ri$ by
$\psi_n(x) = \chi_n(x_1 - \tau(x'))$, where $x = (x_1,x') \in \Ri \times \Ri^{d-1}$ and $\chi_n \colon \Ri \to \Ri$ is 
defined by
\[
\chi_n(t) = \left\{ \begin{array}{ll}
    0 & \mbox{if } t \leq 1/n , \\[5pt]
    \displaystyle   \log( tn)/\log n & \mbox{if } 1/n < t < 1 ,  \\[5pt]
    1 & \mbox{if } t \geq 1 .
            \end{array} \right.
\]
Then $\lim ( \psi_n\varphi) = \one_\Omega \varphi$ in $L_2(\Ri^d)$.
Thus to establish that  $\one_\Omega \varphi \in D(h)$ it suffices to prove  that $ \{ h(\psi_n\varphi) : n \in \Ni \} $ is 
bounded. 
But
\begin{eqnarray*}
h(\psi_n\varphi)
& \leq & 2\, h(\varphi) 
   + 2 \int |\varphi|^2 \sum^d_{k,l=1} c_{kl} \, (\partial_k \psi_n) \, (\partial_l \psi_n)  \\
& \leq & 2\, h(\varphi) 
   + 2 \,(\log n)^{-2}\int_{\Ri^{d-1}} dx'
       \int_{\tau(x') + 1/n}^{\tau(x') + 1} dx_1 \, |\varphi(x)|^2 \, 
     \frac{\langle \nu_x , C(x) \nu_x \rangle}
          {(x_1 - \tau(x'))^2 }
\end{eqnarray*}
for all $n \in \Ni$ where $\nu_x=(1,-(\nabla\tau)(x'))$.
Since the coefficients $c_{kl}$ are in $W^{1,\infty}(\Ri^d)$ there exists an $M > 0$ 
such that 
$|\langle \xi, C(x) \xi \rangle - \langle \xi, C(z) \xi \rangle|
   \leq M \, \|\xi\|^2$
for all $x,z,\xi \in \Ri^d$.
If $x = (x_1,x') \in B_y(r)$, the function $\tau$ is differentiable at $x'$ 
and $x_1 = \tau(x')$ then
\[
\langle \nu_x, C(\tau(x'),x') \nu_x\rangle
   = (1+|(\nabla\tau)(x')|^2)\,\langle n_x, C(\tau(x'),x') n_x \rangle = 0
\]
 by the zero flux condition.
Hence 
$\langle \nu_x,  C(x_1,x') \nu_x \rangle
   \leq M_1 \, |x_1 - \tau(x')|$
for all $x = (x_1,x') \in B_y(r)$ with $\tau$ differentiable at $x'$,
where $M_1 = M (1 + \|\nabla \tau\|_\infty)^2$.
It follows that 
\begin{eqnarray*}
\lefteqn{
(\log n)^{-2}\int_{\Ri^{d-1}} dx'
       \int_{\tau(x') + 1/n}^{\tau(x') + 1} dx_1 \, |\varphi(x_1,x')|^2 \,
     \frac{\langle \nu_x, C(x_1,x') \nu_x\rangle}
          {(x_1 - \tau(x'))^2 }
} \hspace{20mm} \\
& \leq &M_1\,(\log n)^{-2}
       \int_{\Ri^{d-1}} dx'
       \int_{\tau(x') + 1/n}^{\tau(x') + 1} dx_1 \, 
     \frac{|\varphi(x_1,x')|^2 }{(x_1 - \tau(x'))}  
\leq M_1 \,(\log n)^{-1}\|\varphi\|_\infty^2\, |K'|
\end{eqnarray*}
uniformly for all $n \in \Ni$, where $K' \subset \Ri^{d-1}$ is a 
compact set such that $\supp\varphi \subset \Ri \times K'$.
So $ \{ h(\psi_n\varphi) : n \in \Ni \} $ is bounded, as required.
In fact a slightly more detailed argument establishes that 
$\lim h(\psi_n\varphi-\one_\Omega\varphi) = 0$.

\section{Core properties} \label{Sflow4}

In  this section we examine conditions which ensure that $C_c^\infty(\Ri^d)$ is a 
core for the  degenerate elliptic operator $H$ with coefficients $(c_{kl})$ in $W^{1,\infty}$.
Obviously $C_c^\infty(\Ri^d)$ is a core for $H$ if and only if 
$W^{2,\infty}(\Ri^d)$ is a core for $H$.

First, we recall two known core criteria.

\begin{thm} \label{tflow906}
If one of the following two conditions is valid then
$C_c^\infty(\Ri^d)$ is a core for~$H$:
\begin{tabel}
\item \label{tflow906-1}
$c_{kl} \in W^{2,\infty}(\Ri^d)$ for all $k,l \in \{ 1,\ldots,d \} $,
\item \label{tflow906-2}
the matrix $(c_{kl}(x))$ is invertible for all $x \in \Ri^d$.
\end{tabel}
\end{thm}
\proof\
If Condition~\ref{tflow906-1} is valid then $C_c^\infty(\Ri^d)$ is a core
by \cite{Rob7} Section~6, or \cite{ER27} Proposition~2.3, or by 
an adaption of the proof of Proposition~\ref{pord901}.
If Condition~\ref{tflow906-2} is valid then $C_c^\infty(\Ri^d)$ is a core
by the arguments in \cite{Dav14} Theorem~3.1.
Davies requires that the coefficients are smooth, but if the coefficients are bounded
the smoothness condition can be relaxed to
$W^{1,\infty}$.\hfill$\Box$

\ruimte

We shall prove a core theorem with a mixture of the two conditions
of Theorem~\ref{tflow906} in Corollary~\ref{cflow904}.

\begin{lemma} \label{lflow901}
If $\chi \in W^{2,\infty}(\Ri^d)$ and $\varphi \in D(H)$ then 
$\chi \varphi \in D(H)$.
\end{lemma}
Fix $\chi \in W^{2,\infty}(\Ri^d)$.
Then it follows from Lemma~3.4 in \cite{ERSZ2} that $\chi \varphi \in D(h)$ and 
$h(\chi \, \varphi)^{1/2} 
   \leq \|\chi\|_\infty \, h(\varphi)^{1/2} + \|\Gamma(\chi)\|_\infty^{1/2} \, \|\varphi\|_2$
for all $\varphi \in D(h)$, where we define
$\Gamma(\chi) 
= \sum_{k,l=1}^d c_{kl} \, (\partial_k \chi) \, (\partial_l \chi) \in L_\infty$.
If $\varphi,\psi \in C_c^\infty$ then 
\[
h(\psi, \chi \, \varphi)
= h(\chi \, \psi, \varphi) 
   - \sum_{k,l=1}^d \int \psi \, \varphi \, (\partial_k \, c_{kl} \, \partial_l \chi)
   - 2 \sum_{k,l=1}^d \int c_{kl} \, (\partial_k \varphi) \, (\partial_l \chi) \, \psi
\;\;\; .  \]
So 
\begin{equation}
|h(\psi, \chi \, \varphi)|
\leq |h(\chi \, \psi, \varphi)|
   + a \, \|\psi\|_2 \, \|\varphi\|_2 
   + 2 h(\varphi)^{1/2} \, \|\Gamma(\chi)\|_\infty^{1/2} \|\psi\|_2
\;\;\; ,
\label{elflow901;1}
\end{equation}
where $a = \|\sum \partial_k  c_{kl} \partial_l \chi\|_\infty$.
Then by continuity (\ref{elflow901;1}) is valid for all 
$\psi,\varphi \in D(h)$.
Finally, if $\varphi \in D(H)$ then 
$|h(\chi \, \psi, \varphi)| = |(\chi \, \psi, H \varphi)| 
    \leq \|H \varphi\|_2 \, \|\chi\|_\infty \, \|\psi\|_2$
for all $\psi \in D(h)$.
Using (\ref{elflow901;1}) it follows that there exists a $c > 0$ 
such that $|h(\psi, \chi \varphi)| \leq c \, \|\psi\|_2$ 
for all $\psi \in D(h)$.
Therefore $\chi \varphi \in D(H)$.\hfill$\Box$

\ruimte

If $A \subset \Ri^d$ with $A \neq \emptyset$ and $\delta > 0$ 
define the open set $A_\delta \subset \Ri^d$ by
$A_\delta = \{ x \in \Ri^d : d(x,A) < \delta \} $.

\begin{lemma} \label{lflow902}
Let $H_1$ and $H_2$ be degenerate elliptic operators with $W^{1,\infty}$-coefficients 
$(c^{(1)}_{kl})$ and $(c^{(2)}_{kl})$ and let $h^{(1)}$ and $h^{(2)}$
be the corresponding quadratic forms.
Let $U \subset \Ri^d$ be an open set and suppose that 
$c^{(1)}_{kl}|_U = c^{(2)}_{kl}|_U$ for all $k,l \in \{ 1,\ldots,d \} $.
Let $\varphi \in L_2(\Ri^d) \setminus \{ 0 \} $ 
and suppose that $(\supp \varphi)_\delta \subset U$.

Then $\varphi \in D(h^{(1)})$ if and only if $\varphi \in D(h^{(2)})$
and then $h^{(1)}(\varphi) = h^{(2)}(\varphi)$.
Similarly, $\varphi \in D(H_1)$ if and only if $\varphi \in D(H_2)$
and then $H_1 \varphi = H_2\varphi$.
Moreover, $\supp H_1 \varphi \subseteq \supp \varphi$.
\end{lemma}
\proof\
There exists a $\chi \in W^{2,\infty}(\Ri^d)$ such that 
$\chi|_{\supp \varphi} = \one$ and $\supp \chi \subset U$.
Suppose $\varphi \in D(h^{(1)})$.
Then there exists a sequence $\varphi_1,\varphi_2,\ldots \in W^{1,2}(\Ri^d)$
such that $\lim \varphi_n = \varphi$ in $D(h^{(1)})$.
Then $\lim \varphi_n = \varphi$ in $L_2(\Ri^d)$.
But $h^{(1)}(\chi \varphi_n) = h^{(2)}(\chi \varphi_n)$ and 
$h^{(1)}(\chi \varphi_n - \chi \varphi_m) = h^{(2)}(\chi \varphi_n - \chi \varphi_m)$
for all $n,m \in \Ni$.
Therefore $\chi \varphi_1,\chi \varphi_2$ is a Cauchy sequence in $D(h^{(2)})$.
Since $\lim \chi \varphi_n = \varphi$ in $L_2$ one deduces that 
$\varphi \in D(h^{(2)})$ and 
$h^{(2)}(\varphi) = h^{(1)}(\varphi)$.

Finally suppose that $\varphi \in D(H_1)$.
If $\psi \in C_c^\infty(\Ri^d)$ with $\supp \psi \subset (\supp \varphi)^{\rm c}$
then $(H_1 \varphi,\psi) = h^{(1)}(\varphi,\psi) = 0$ by locality.
Therefore $\supp H_1 \varphi \subseteq \supp \varphi$.
Clearly $\varphi \in D(h^{(1)})$ and by the first part, also 
$\varphi \in D(h^{(2)})$.
Let $\psi \in D(h^{(2)})$.
Then $\chi \psi \in D(h^{(2)})$ and $\supp \chi \psi \subset U$.
Therefore $\chi \psi \in D(h^{(1)})$.
Then by locality one deduces that 
$h^{(2)}(\varphi,\psi) 
   = h^{(2)}(\varphi, \chi \psi) + h^{(2)}(\varphi, (\one - \chi) \psi)
   = h^{(2)}(\varphi, \chi \psi)
   = h^{(1)}(\varphi, \chi \psi)$.
So $|h^{(2)}(\varphi,\psi)| = |h^{(1)}(\varphi, \chi \psi)| = |(H_1 \varphi, \chi \psi)|
   \leq \|H_1 \varphi\|_2 \, \|\chi\|_\infty \, \|\psi\|_2$.
Therefore $\varphi \in D(H_2)$.
If $\psi \in C_c^\infty(U)$ then 
$(H_1 \varphi, \psi) 
   = (\varphi, H_1 \psi)
   = (\varphi, H_2 \psi)
   = (H_2 \varphi, \psi)$.
Since $\supp H_1 \varphi \subseteq U$ and $\supp H_2 \varphi \subseteq U$
it follows that $H_1 \varphi = H_2 \varphi$.\hfill$\Box$

\begin{prop} \label{pflow903}
Let $A \subset \Ri^d$, $\delta > 0$, 
let $H_1$ and $H_2$ be degenerate elliptic operators with 
$W^{1,\infty}$-coefficients 
$(c^{(1)}_{kl})$ and $(c^{(2)}_{kl})$.
Suppose $\emptyset \neq A \neq \Ri^d$, 
$c^{(1)}_{kl}|_{A_\delta} = c_{kl}|_{A_\delta}$ and
$c^{(2)}_{kl}|_{(A^{\rm c})_\delta} = c_{kl}|_{(A^{\rm c})_\delta}$ 
for all $k,l \in \{ 1,\ldots,d \} $ and 
$C_c^\infty(\Ri^d)$ is a core for both $H_1$ and $H_2$.
Then $C_c^\infty(\Ri^d)$ is a core for $H$.
\end{prop}
\proof\
Let $\tau \in C_c^\infty(\Ri^d)$ be such that $\int \tau = 1$ and 
$\tau(x) = 0$ for all $x \in \Ri^d$ with $|x| > \frac{\delta}{4}$.
Let $\chi = \tau * \one_{A_{\delta / 2}}$.
Then $\chi \in W^{2,\infty}(\Ri^d)$,
$\chi|_{A_{\delta / 4}} = \one$ and $\supp \chi \subset A_{3 \delta / 4}$.
Moreover, $\supp (\one - \chi) \subset (A_{\delta / 4})^{\rm c} \subset A^{\rm c}$.
There exist $\chi_1,\chi_2 \in W^{\infty,\infty}(\Ri^d)$
such that 
$\chi_1|_{A_{3 \delta / 4}} = \one$, $\supp \chi_1 \subset A_\delta$,
$\chi_2|_{A^{\rm c}} = \one$ and $\supp \chi_2 \subset (A^{\rm c})_\delta$.

Let $\varphi \in D(H)$.
It follows from Lemma~\ref{lflow901} that $\chi \varphi \in D(H)$ and 
$(\one - \chi) \varphi \in D(H)$.
We shall show that we can approximate both elements by $C_c^\infty$-functions.
We may assume that $\chi \varphi \neq 0 \neq (\one - \chi) \varphi$.
Since $\supp (\chi \varphi) \subset A_{3 \delta / 4}$ one deduces 
from Lemma~\ref{lflow902} that 
$\chi \varphi \in D(H_1)$ and $H_1(\chi \varphi) = H(\chi \varphi)$.
By assumption there exist $\varphi_1,\varphi_2,\ldots \in C_c^\infty(\Ri^d)$
such that 
$\lim \varphi_n = \chi \varphi$ in $D(H_1)$.
Then $\lim \chi_1 \varphi_n = \chi_1 \chi \varphi = \chi \varphi$ in $D(H_1)$ by 
Lemma~\ref{lflow901}.
But $\chi_1 \varphi_n \in C_c^\infty(\Ri^d)$ and
$\supp \chi_1 \varphi_n \subset A_\delta$ for all $n \in \Ni$.
Therefore $\chi_1 \varphi_n \in D(H)$ and 
$H(\chi_1 \varphi_n) = H_1(\chi_1 \varphi_n)$, again by 
Lemma~\ref{lflow902}.
So $\lim \chi_1 \varphi_n = \chi \varphi$ in $D(H)$.
Similarly, using $H_2$ and $\chi_2$ there exists a sequence
$\psi_1,\psi_2,\ldots \in C_c^\infty(\Ri^d)$ such that 
$\lim \chi_2 \psi_n = (\one - \chi) \varphi$ in $D(H)$.
Then $\lim (\chi_1 \varphi_n + \chi_2 \psi_n) = \varphi$ in $D(H)$.
Since $\chi_1 \varphi_n + \chi_2 \psi_n \in C_c^\infty(\Ri^d)$ the
proposition follows.\hfill$\Box$

\begin{cor} \label{cflow904}
Suppose there exist a set $A$ and $\delta > 0$ such that 
$\emptyset \neq A \neq \Ri^d$, the 
matrix $(c_{kl}(x))$ is invertible for all $x \in (A^{\rm c})_\delta$
and $c_{kl}|_{A_\delta} \in W^{2;\infty}(A_\delta)$.
Then $C_c^\infty(\Ri^d)$ is a core for~$H$.
\end{cor}
\proof\
There exists a $\chi_1 \in W^{2,\infty}(\Ri^d)$ such that 
$\chi_1|_{A_{\delta/2}} = \one$ and $\supp \chi_1 \subset A_\delta$.
Define $c^{(1)}_{kl} = \chi_1 \, c_{kl} \in W^{2,\infty}(\Ri^d)$.
Then $c^{(1)}_{kl}|_{A_{\delta/2}} = c_{kl}|_{A_{\delta/2}}$.

There exists a $\chi_2 \in W^{1,\infty}(\Ri^d)$ such that 
$\chi_2|_{(A^{\rm c})_{\delta/2}} = \one$ and $\supp \chi_2 \subset (A^{\rm c})_\delta$.
Define 
$c^{(2)}_{kl} = \chi_2 \, c_{kl} + (\one - \chi_2) \delta_{kl} \in W^{1,\infty}(\Ri^d)$.
Let $H_1$ and $H_2$ be the degenerate elliptic operator with coefficients $(c_{kl}^{(1)})$
and $(c_{kl}^{(2)})$.
Now apply Theorem~\ref{tflow906}.\ref{tflow906-1} to $H_1$,
Theorem~\ref{tflow906}.\ref{tflow906-2} to $H_2$ and 
use Proposition~\ref{pflow903}.\hfill$\Box$

\subsection*{Acknowledgement}
Part of this  work was carried out whilst the first author was visiting 
the Australian National University
with partial support from the Centre for Mathematics and its
Applications and part of the work was
carried out whilst the second author was visiting the University of Auckland 
with financial support
from the Faculty of Science.

\end{document}